\documentclass[11pt,a4paper]{amsart}
\usepackage{graphicx,multirow,array,amsmath,amssymb,kotex}

\usepackage{xcolor}
\usepackage{longtable}
\usepackage{supertabular} 

\newtheorem{theorem}{Theorem}

\newtheorem*{definition}{Definition}

\begin{document}
	
\title{Folded ribbonlength of 2-bridge knots}
	
\author[H. Kim]{Hyoungjun Kim}
\address{College of General Education, Kookmin University, Seoul 02707, Korea}
\email{kimhjun@kookmin.ac.kr}
\author[S. No]{Sungjong No}
\address{Department of Mathematics, Kyonggi University, Suwon 16227, Korea}
\email{sungjongno@kgu.ac.kr}
\author[H. Yoo]{Hyungkee Yoo}
\address{Research Institute for Natural Science, Hanyang University, Seoul 04763, Korea}
\email{hyungkee@hanyang.ac.kr}

	
\begin{abstract}
A ribbon is a two-dimensional object with one-dimensional properties which is related with geometry, robotics and molecular biology.
A folded ribbon structure provides a complex structure through a series of folds.
We focus on a folded ribbon with knotted core.
The folded ribbonlength $Rib(K)$ of a knot $K$ is the infimum of the quotient of length by width among the ribbons representing a knot type of $K$.
This quantity tells how efficiently the folded ribbon is realized.
Kusner conjectured that folded ribbonlength is bounded by a linear function of the minimal crossing number $c(K)$.
In this paper, we confirm that the folded ribbonlength of a 2-bridge knot $K$ is bounded above by $2c(K)+2$.
\end{abstract}

\maketitle

\section{Introduction} \label{sec:intro}

Origami, a Japanese word meaning paper-folding, is the process of making a two-dimensional object into a three-dimensional object.
Recently, it has been attracting attention in the fields of geometry, robotics, and molecular biology over the past two decades~\cite{BM,CTZ,DM,HABTKDW,SALTSA}.
A ribbon is a two-dimensional object with one-dimensional properties.
A folded ribbon structure provides a complex structure through a series of folds~\cite{CDBG,RLRA}.
For example, a folded ribbon realizes a ribosomal walking robot~\cite{RCL,WPF}.
DNA is the famous ribbon shaped object found in nature.
A knotted ribbon shape is observed in many circular DNA~\cite{LJ,LLA,LPCW,SCBLC}.
Recently, there is a result that it is possible to create a wider ribbon structure by connecting multiple DNA in parallel~\cite{HPLY}.
We focus on representing the knotted ribbon in a folded form.

\begin{figure}[h!]
    \includegraphics{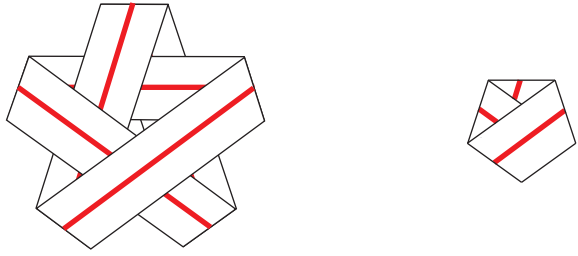}
    \caption{Folded ribbon knots of the trefoil knot}
    \label{fig:frk}
\end{figure}

A {\it knot} is a simple closed curve in $\mathbb{R}^3$, and a {\it link} is a disjoint union of knots.
For convenience, we refer to the knot and link together as a knot in the remaining part of this paper.
Kauffman~\cite{K} introduced a {\it folded ribbon knot}, which is a representation of a knot that folds a long rectangular strip to become flat.
Because of the construction of this representation, the core of the ribbon is piecewise linear as drawn in Figure~\ref{fig:frk}.
Since the right figure of Figure~\ref{fig:frk} is represented tighter than the left figure, the right figure can be represented with a shorter ribbon than the left figure.
Moreover, the length of folded ribbon knot depends on the width of the ribbon.
Thus we consider the ratio of the length and width of the ribbon to see how efficiently the knotted ribbon is represented.

Now we define a folded ribbon knot and a folded ribbonlength which are based on~\cite{DKTZ}.

\begin{definition}
Let $K$ be the knot. 
Then $K_w$ is a folded ribbon knot of $K$  with width $w$ that satisfies following two conditions:
\begin{enumerate}
    
\item The ribbon is flat and its fold lines are disjoint.
\item The core of $K_w$ is a union of a finite number of consecutively straight lines with crossing information which represents the knot type of $K$.

\end{enumerate}

\end{definition}

Note that the trivial knot can be represented as a folded ribbon obtained by folding twice so that these folding lines are perpendicular to the core.
In this case, the length of the ribbon can be made arbitrarily small.
So we use infimum when defining the folded ribbonlength.

\begin{definition}
Let $K$ be a knot (or link) and $w$ be a positive real number.
Then
$$Rib(K)=\inf_{K' \in [K]_w}\frac{\emph{\text{Len}}(K')}{w}$$
is called a folded ribbonlength of $K$
where $[K]_w$ is the set of folded ribbon knots with width $w$ representing a knot type $K$, and ${\emph{\text{Len}}(K')}$ is the length of the core of $K'$.
\end{definition}

In 2005, Kauffman~\cite{K} introduced models for folded ribbon knots of the trefoil and the figure eight knot, and their folded ribbonlength.
In the paper, Kusner's conjecture is reported, which proposes a linear relation between the folded ribbonlength and the crossing number of a knot $K$:
$$c_1 \cdot c(K) \le Rib(K) \le c_2 \cdot c(K)$$
In 2008, Kennedy et al~\cite{KMRT} provided bounds for the constants $c_1$ and $c_2$ by calculating the folded ribbonlength of torus knots.

Many researchers have found upper bounds for the folded ribbonlength of knots. In 2017, Tian~\cite{T} found the following upper bound 
$$Rib(K) \leq 2c(K)^2 + 6c(K) + 4,$$
and in 2020, Denne~\cite{D} improved the result.
Furthermore, Denne et al.~\cite{DHLM} provided linear upper bounds for the folded ribonlength of 2-bridge knots, $(2,p)$-torus knots and $(p,q,r)$-pretzel knots in terms of their crossing numbers. 
In particular their upper bound for 2-bridge knots is as follows. 
$$Rib(K) \leq 6c(K)-2.$$

In this paper we improve the result of the folded ribbonlength of 2-bridge knots.

\begin{theorem}
Let $K$ be a 2-bridge knot or link.
Then the upper bound of the folded ribbonlength of $K$ is $2c(K)+2$.
\end{theorem}

\section{2-bridge knots} \label{sec:bvl}

\begin{figure}[h!]
    \includegraphics[scale=0.9]{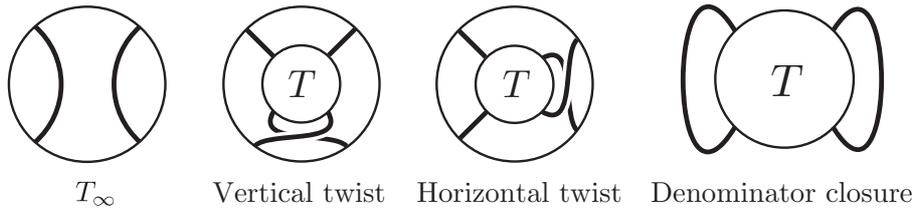}
    \caption{Operations for rational tangles}
    \label{fig:tangle}
\end{figure}

Conway~\cite{C} defined a 2-tangle as a portion of knot diagram from which there emerge just 4 arcs pointing in the compass directions NW, NE, SW and SE.
Figure~\ref{fig:tangle} shows $T_{\infty}$, a vertical twist, a horizontal twist and a denominator closure.
The tangle $T_{\infty}$ is one of the simplest 2-tangles which consists of two vertical strings.
A {\it vertical twist} is an operation of a 2-tangle which switches the location of two endpoints SW and SE by using an rotation around N-S axis counterclockwise.
A {\it horizontal twist} is an operation of a 2-tangle which switches the location of two endpoints NE and SE by using an rotation around E-W axis clockwise.
Especially, a tangle which only consists of vertical twists (or horizontal twists) is called an {\it integer tangle}.
A 2-tangle that is obtained from $T_{\infty}$ by a sequence of vertical and horizontal twists and their inverses is called a {\it rational tangle}.
We can describe a rational tangle by using a such sequence.
A rational tangle with Conway notation $T(a_1,\dots,a_m)$ means start with $T_{\infty}$, then using $a_1$ vertical twists, then $a_2$ horizontal twists, and repeat until $a_m$.
The fraction of a rational $q$-tangle $T(a_1,\dots,a_m)$ is then defined as the number given by the continued fraction 
$$q=a_m + \cfrac{1}{a_{m-1} + \cfrac{1}{ \ddots +\cfrac{1}{a_2+ \cfrac{1}{a_1} }}}$$
where $q$ is a rational number.
Conway~\cite{C} showed that two rational tangles are equivalent if and only if the corresponding continued fractions have the same value.
In~\cite{BZ}, it is shown that for any rational $q$-tangle, there is a corresponding Conway notation which consists of an odd number of positive integers (or negative integers) when $q$ is positive (or negative).

\begin{figure}[h!]
    \includegraphics[scale=0.8]{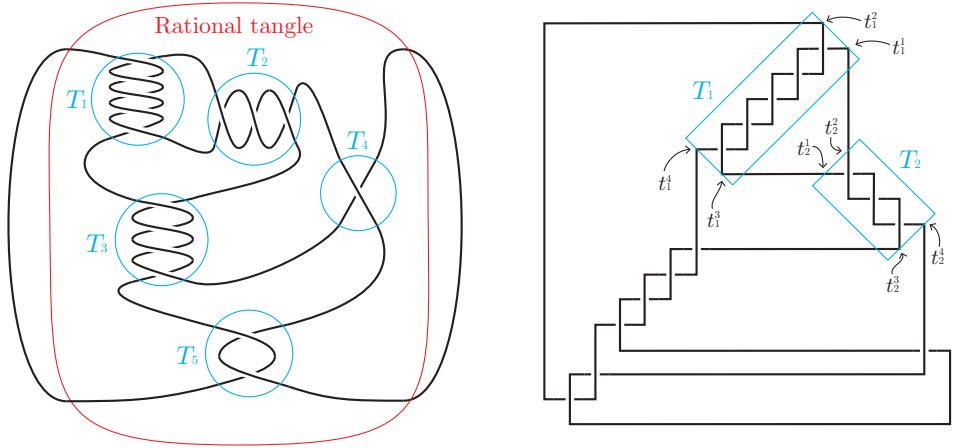}
    \caption{Two diagrams of a 2-bridge knot obtained from a rational tangle $T(-5,-3,-4,-1,-2)$}
    \label{fig:2br}
\end{figure}

A 2-bridge knot is also called a {\it rational knot}, that is obtained from a rational tangle by taking a denominator closure which connects NE to SE, and NW to SW.
Figure~\ref{fig:2br} is an example of a 2-bridge knot which is obtained from a rational tangle $T(-5,-3,-4,-1,-2)$.
Since the right figure is a suitable diagram to show the process of the construction of a rational tangle, we construct a folded ribbon knot based on this diagram.
Note that a diagram of a rational knot which is obtained from a rational tangle whose Conway notation consists of an odd number of positive integers (or negative integers) is a non-nugatory alternating diagram.
Thus this diagram has the minimum number of crossings by~\cite{K1,Mu,Th}.

To simplify notation in the proof, we define new notation related to rational tangles.
Let $T_i$ denote the $i$-th integer tangle in a rational tangle that corresponds to $a_i$ in the Conway notation as drawn in Figure~\ref{fig:2br}.
Each integer tangle $T_i$ has four endpoints, and notations of these endpoints are differently defined depend on the value $i$, since the shape of $T_i$ is different depending on whether $i$ is odd or even.
If $i$ is odd, then $t_i^1$, $t_i^2$, $t_i^3$ and $t_i^4$ are four endpoints NE, NW, SE and SW of $T_i$, respectively.
If $i$ is even, then $t_i^1$, $t_i^2$, $t_i^3$ and $t_i^4$ are four endpoints SW, NW, SE and NE of $T_i$, respectively.
In the right figure of Figure~\ref{fig:2br}, we easily check that four endpoints of $T_i$ are labeled in order from right to left when $i$ is odd, and from left to right when $i$ is even.
When $i$ is odd, the shape of $T_i$ is symmetrical to the shape when $i$ is even.
Thus, the definitions above are convenient for construction of a folded ribbon knot.

Note that a rational knot whose Conway notation consists of an odd number of negative integers has an alternating diagram with the minimum number of crossings provided its rational tangle has a negative value. For the figures used in the proof, it’s easy to recognize when all the numbers in the Conway notation are negative.
So, we will focus on the construction of a folded ribbon knot whose Conway notation consists of negative integers.
In the case of a rational knot whose rational tangle has positive value, we construct its folded ribbon knot in a precisely symmetrical way.

\section{Construction of 2-bridge folded ribbon knot} \label{sec:main}

\begin{figure}[h!]
    \includegraphics{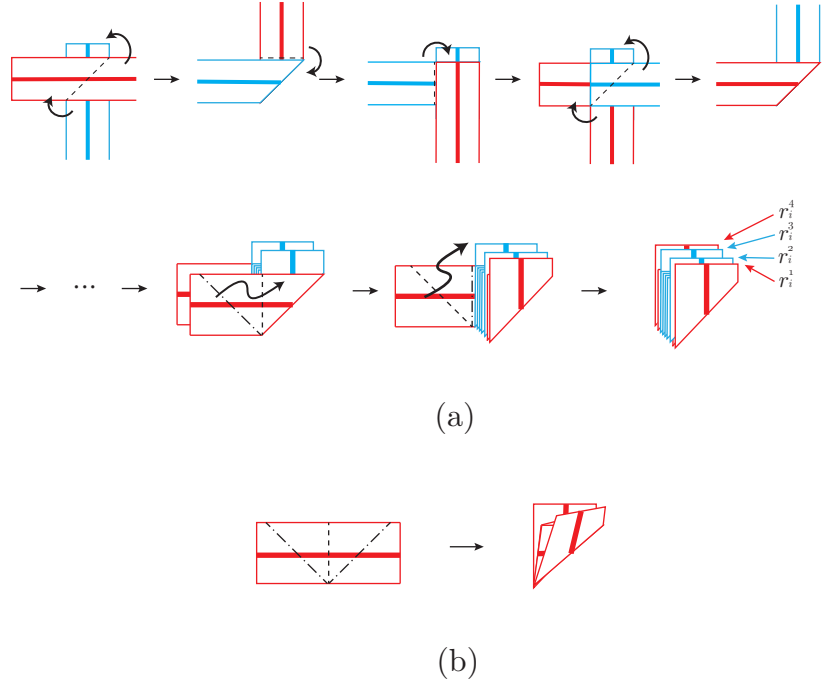}
    \caption{Construction of integer ribbons}
    \label{fig:integer1}
\end{figure}

In this section, we construct a 2-bridge folded ribbon knot step by step.
As mentioned in the previous section, we construct a 2-bridge knot whose Conway notation consists of negative integers.
To construct a 2-bridge knot, first we construct an integer tangle $T_i$ by folded ribbons, which is called the {\it integer ribbon} $R_i$ as drawn in Figure~\ref{fig:integer1}~(a).
Assume that $T_i$ consists of $-n$ half twists for a positive integer $n$.
We use two ribbons with width $w$ such that lengths of the two ribbons are $nw+2 \epsilon$ when $n$ is even, and lengths of the two ribbons are $(n+1)w+2 \epsilon$ and $(n-1)w+2 \epsilon$ respectively when $n$ is odd, for small $\epsilon > 0$.
This means that the sum of the lengths of the two ribbons is $2nw+4\epsilon$.
Here, the actual length of the ribbon used to construct the integer ribbon is $2nw$, but for visibility of the connecting process of integer ribbon, leave a margin of length $\epsilon$ at each end of the ribbons.
First, place one ribbon vertically on the plane.
Here, if the two ribbons are of different lengths, place the shorter one.
Next, place the other ribbon horizontally such that the length of a portion of the horizontal ribbon located to the left of the vertical ribbon is $w+\epsilon$, and the length of a portion of the vertical ribbon located above the horizontal ribbon is $\epsilon$.
Then, fold the two ribbons together from the bottom right to the top left as drawn in the first part of Figure~\ref{fig:integer1}~(a).
During the process as depicted in first row in Figure~\ref{fig:integer1}~(a), the positions of the two ribbons are changed and the length of each ribbon is shortened by $w$.
Repeat this until the remaining length of one ribbon is $\epsilon$.
Finally fold up the remaining parts of the ribbon as drawn in the last row of Figure~\ref{fig:integer1}~(a).
Remark that the resulting ribbon is embedded in an isosceles right triangle except for the four end lines of two ribbons.
We label for these end lines as $r_i^1,\dots,r_i^4$ in order from right to left.

Note that if $n$ is 1, then the integer ribbon is constructed by one ribbon with length $2w+2 \epsilon $ as drawn in Figure~\ref{fig:integer1}~(b).
For the case that $n$ is negative integer, consider the mirror symmetry of the above process.
Then four end lines $r_i^1$, $r_i^2$, $r_i^3$ and $r_i^4$ correspond to the four endpoints $t_i^1$, $t_i^2$, $t_i^3$ and $t_i^4$, respectively.

\begin{figure}[h!]
    \includegraphics{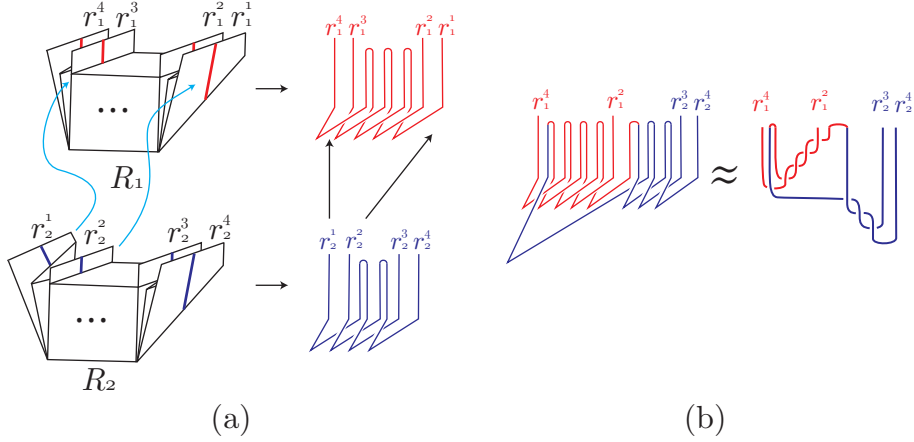}
    \caption{Connecting two integer ribbons $R_1$ and $R_2$}
    \label{fig:integer2}
\end{figure}

Now we construct the folded ribbon knot of a 2-bridge knot $K$ by connecting integer ribbons.
First, connect $R_1$ and $R_2$ by attaching $r_1^1$ to $r_2^2$, and $r_1^3$ to $r_2^1$.
To attach $r_1^1$ to $r_2^2$, we put the part of $R_2$ containing $r_2^1$ between $r_1^3$ and $r_1^4$, and the remaining part $R_2$ to the right of $r_1^1$ as drawn in Figure~\ref{fig:integer2}~(a).
Even though the integer ribbon $R_1$ has $|a_1|-1$ crossings, one crossing is added to $R_1$ in the connecting process for $R_1$ and $R_2$ as drawn in Figure~\ref{fig:integer2} (b).
Especially, if $a_1=-1$, then we use the labeling $r_1^2$ instead of $r_2^1$ after connecting $R_1$ and $R_2$.
Similarly, if $a_2=-1$, then we use the labeling $r_2^3$
instead of $r_1^1$.
In Figure~\ref{fig:integer2} (a), the left figures show integer ribbons, and the right figures show their cores.
Figure~\ref{fig:integer2} (b) shows the core of the result after connecting two integer ribbons.

\begin{figure}[h!]
    \includegraphics{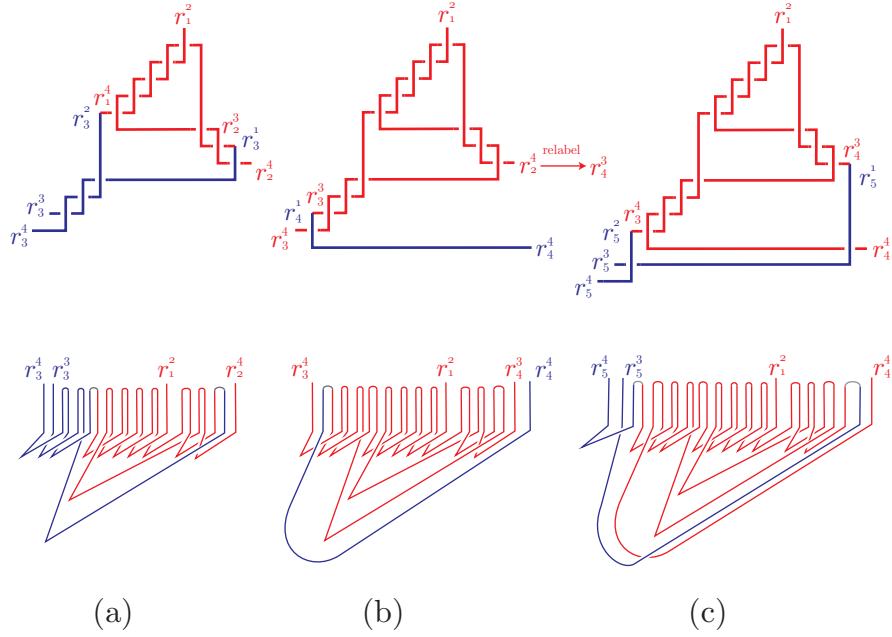}
    \caption{Intermediate processes of consecutive connections of integer ribbons}
    \label{fig:integer3}
\end{figure}

We consecutively connect the remaining integer ribbons $R_i$ for $i \geq 3$ in the similar way.
In detail, to attach $r_{i-1}^3$ and $r_i^1$, we put the part of $R_i$ containing $r_i^1$ between $r_{i-1}^3$ and $r_{i-1}^3$.
Moreover to attach $r_{i-2}^4$ to $r_i^2$, we put the remaining part $R_i$ to the left of $r_{i-2}^4$ when $i$ is odd, and to the right otherwise as drawn in Figure~\ref{fig:integer3}.
Especially, if $a_i=-1$, then we use the labeling $r_i^3$ instead of $r_{i-2}^4$ after connecting $R_i$ as drawn in Figure~\ref{fig:integer3} (b) and (c).

\begin{figure}[h!]
    \includegraphics{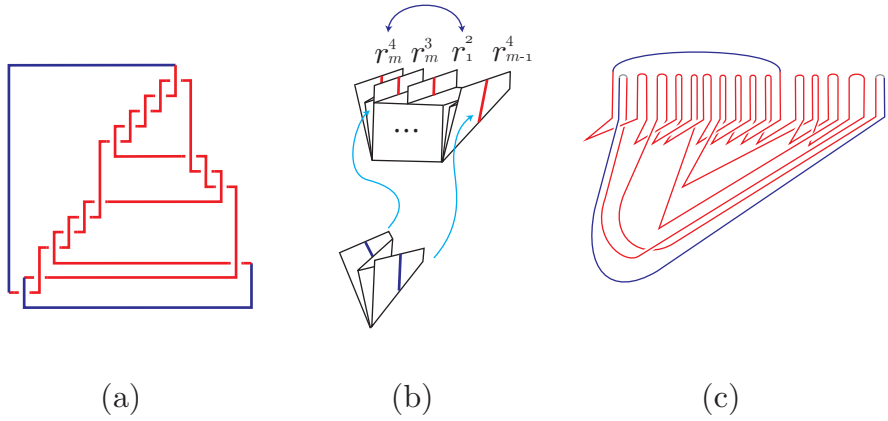}
    \caption{The final process to construct a folded ribbon 2-bridge knot}
    \label{fig:integer4}
\end{figure}

So far we have constructed a rational tangle of $K$ in the form of a folded ribbon.
Now all that remains is taking a denominator closure for this tangle to complete the 2-bridge knot $K$ as drawn in Figure~\ref{fig:integer4} (a).
Note that the denominator closure consists of two arcs, one connecting $t_1^2$ and $t_m^4$, and the other connecting $t_{m-1}^4$ and $t_m^3$.
First connect $r_{m-1}^4$ and $r_m^3$ by attaching a new integer ribbon with length $2w$.
Then we can attach $r_1^2$ to $r_m^4$ directly as drawn in Figure~\ref{fig:integer4} (b).

Finally, find the length of folded ribbon required to construct $K$.
Note that the length of integer tangle $R_i$ is at most $2nw+4\epsilon$ when its corresponding tangle $T_i$ has $n$ half twists.
So the length of each integer ribbon $R_1,\dots,R_m$ is bounded above the number of half twists of the corresponding integer tangles $T_1,\dots,T_m$ multiplied by $2w$ and then added by $4\epsilon$.
Since the sum of the half twists of all integer tangles $T_1,\dots,T_m$  is equal to the crossing number $c(K)$ of $K$, the length of the folded ribbon representing the rational tangle is bounded above $2wc(K)+4m\epsilon$.
Furthermore, the denominator closure in the form of a folded ribbon has length $2w+2\epsilon$.
Thus the length of folded ribbon knot is equal to $2(c(K)+1)w+(4m+2)\epsilon$.
Therefore, the upper bound of a folded ribbonlength of a 2-bridge knot is $2(c(K)+1)$ when $\epsilon$ goes to 0.

\section*{acknowledgement}

The first author(Hyoungjun Kim) was supported by the National Research Foundation of Korea (NRF) grant funded by the Korea government Ministry of Science and ICT(NRF-2021R1C1C1012299).
The corresponding author(Sungjong No) was supported by the National Research Foundation of Korea(NRF) grant funded by the Korea government Ministry of Science and ICT(NRF-2020R1G1A1A01101724).

\end{document}